\documentclass[10pt]{article}
\usepackage{amssymb, amsmath, url, graphicx, setspace, geometry}

\def\3{\subset }
\def\4{\subseteq }
\def\<{\left<}
\def\>{\right>}

\def\bit{\begin{itemize}}
\def\eit{\end{itemize}}
\def\3{\subset }
\def\4{\subseteq }

\def\0{\leqno}

\def\barr{\begin{array}}
\def\earr{\end{array}}

\def\Z{{\rlap{$\kern2pt{\rm Z}$}{\rm Z}\,}}
\def\bld#1#2{{\buildrel{#1}\over{#2}}}
\def\st#1#2{{\mathrel{\mathop{#2}\limits_{#1}}{}\!}}
\def\stb#1#2#3{{\st{{#1}}{\bld{{#2}}{#3}}{}\!}}
\def\xmare#1#2{\stb{#1}{#2}{\mbox{\Huge$\times$}}}

\title{\bf A note on the number of cyclic subgroups of a finite group}
\author{Marius T\u arn\u auceanu and Mihai-Silviu Lazorec}
\date{May 1, 2018}

\begin{document}

\maketitle

\begin{abstract}
Let $G$ be a finite group, $L_1(G)$ be its poset of cyclic subgroups and consider the quantity $\alpha(G)=\frac{|L_1(G)|}{|G|}$. The aim of this paper is to study the class $\cal{C}$ of finite nilpotent groups having $\alpha(G)=\frac{3}{4}$. We show that if $G$ belongs to this class, then it is a 2-group satisfying certain conditions. Also, we study the appartenance of some classes of finite groups to $\cal{C}$.
\end{abstract}



\section{Introduction}

Let $G$ be a finite group and $L_1(G)$ be its poset of cyclic subgroups. In \cite{4}, the quantity $$\alpha(G)=\frac{|L_1(G)|}{|G|}$$ was introduced and studied, the authors indicating and proving relevant properties of this number. Some of their results will be used in this paper and we recall the following ones:
\begin{itemize}
\item[--] If $(G_i)_{i=\overline{1,k}}$ is a family of finite groups having coprime orders, then
$$\alpha(\xmare{i=1}k G_i)=\prod\limits_{i=1}^k \alpha(G_i).$$
\item[--] $\alpha(G)=1$ if and only if $G$ is an elementary abelian 2-group.
\item[--] If $N\triangleleft\, G$ then $\alpha(G)\leq \alpha(\frac{G}{N})$, and if the equality holds, then $N$ is an elementary abelian 2-group.
\item[--] $\alpha(G)=\alpha(G\times \mathbb{Z}_2^n), \ \forall \ n\in\mathbb{N}.$
\item[--] The finite groups $G$ having $\alpha(G)>\frac{3}{4}$ were entirely classified (see \textbf{Theorem 5} of \cite{4}). We remark that the elementary abelian 2-groups are the only finite abelian groups that satisfy this property.
\item[--] The value $\frac{3}{4}$ is the largest non-trivial accumulation point of the set $\lbrace \alpha(G) \ | \ G=\text{finite group} \rbrace.$
\end{itemize}
Having in mind especially the last two properties, it is natural to ask if it is possible to fully describe the finite groups having $\alpha(G)=\frac{3}{4}$. A starting point for an answer is given by this paper in which we will study the following class of groups:
$$\mathcal{C}=\lbrace G=\text{finite nilpotent group} \ | \ \alpha(G)=\frac{3}{4}\rbrace.$$

The paper is organized as follows. Section 2 deals with the appartenance of finite abelian groups to $\mathcal{C}$. We prove our main result which states that if a group $G$ is contained in $\mathcal{C}$, then it is a 2-group with $G'=\Phi(G)$ or there is $n\in\mathbb{N}$ such that $\frac{G}{G'}\cong \mathbb{Z}_2^n\times \mathbb{Z}_4$ and $G'$ is elementary abelian. The aim of Section 3 is to investigate if there are any finite groups contained in $\mathcal{C}\cap \mathcal{G}$, where $\mathcal{G}$ is the class of (almost) extraspecial 2-groups, the class of generalized dicyclic 2-groups, the class of generalized dihedral 2-groups  and the class of $2$-groups possessing a cyclic maximal subgroup, respectively. We end our paper by suggesting some open problems in the last Section.

Most of our notation is standard and will usually not be repeated here. Elementary notions and results on groups can be found in \cite{5,9}. For subgroup lattice concepts we refer the reader to \cite{8,10,11}.

\section{A description of the finite groups contained in $\cal{C}$}

Let $m$ be a positive integer. For a finite $p$-group $G$ of exponent $p^m$, we denote by $n_{p^i}(G)$ the number of cyclic subgroups of order $p^i$ of $G$, where $i=1,2,\ldots,m.$ The first step that need to be done before we can provide the features of the groups which belong to $\cal{C}$ is to find the abelian groups that are contained in $\cal{C}$. To complete this task, we will prove two preliminary results. The first one mainly states that the number of cyclic subgroups of a finite $p$-group of order $p^n$ is less than or equal to the number of cyclic subgroups of $\mathbb{Z}_p^n$. The second result establishes an inequality that involves $\alpha(G)$, where $G$ is a finite $p$-group of order $p^n$ and $p$ is an odd prime number.\\

\textbf{Lemma 2.1.} \textit{Let $n$ be a positive integer and let $G$ be a finite $p$-group of order $p^n$. Then $\alpha(G)\leq\alpha(\mathbb{Z}_p^n).$}

\textbf{Proof.} Let $m$ be a positive integer such that $exp(G)=p^m$. Then the order of $G$ is
$$p^n=1+\sum\limits_{i=1}^mn_{p^i}(G)\varphi(p^i).$$
Consequently, we have
$$p^n-1=\sum\limits_{i=1}^m n_{p^i}(G)\varphi(p^i)\ge (n_p(G)+n_{p^2}(G)+\ldots +n_{p^m}(G))(p-1).$$ This implies that
$$\alpha(G)=\frac{1+n_p(G)+\ldots+n_{p^m}(G)}{p^n}\le\frac{1+\frac{p^n-1}{p-1}}{p^n}=\alpha(\mathbb{Z}_p^n),$$ as desired.
\hfill\rule{1,5mm}{1,5mm} \\

\textbf{Lemma 2.2.} \textit{Let $n$ be a positive integer, $p$ be an odd prime number and $G$ be a finite $p$-group of order $p^n$. Then $\alpha(G)<\frac{3}{4}$.}

\textbf{Proof.} For $n=1$, it is obvious that $\alpha(G)=\frac{2}{p}<\frac{3}{4}$ since $p\ge3$. Let $n\ge2$ be a positive integer. According to \textbf{Lemma 2.1}, we have $\alpha(G)\leq\alpha(\mathbb{Z}_p^n)$, so it suffices to show that $\alpha(\mathbb{Z}_p^n)<\frac{3}{4}$. The last inequality is equivalent to $$3p^{n+1}-7p^n-4p+8>0.$$
But, the above relation is true for $p\ge 3$ and  $n\ge 2$ since $$3p^{n+1}-7p^n-4p+8=(3p-7)p^n-4p+8\ge 2p^n-4p+8>0.$$
\hfill\rule{1,5mm}{1,5mm} \\

Let $G$ be a finite abelian group and denote by $L(G)$ its subgroup lattice. Remark that, in this case, the function $\alpha:L(G)\longrightarrow [0,1]$ is decreasing. Let $H$ and $K$ be two subgroups of $G$ such that $H\subseteq K$. Since $G$ is abelian, there is a subgroup $L$ of $K$ such that $H\cong \frac{K}{L}$. Then $\alpha(H)=\alpha(\frac{K}{L})\ge\alpha(K)$. We have all the necessary ingredients to determine the finite abelian groups that are contained in $\mathcal{C}$.\\

\textbf{Theorem 2.3.} \textit{The only finite abelian groups that belong to $\mathcal{C}$ are $\mathbb{Z}_2^n\times\mathbb{Z}_4,$ where $n\in\mathbb{N}$.}

\textbf{Proof.} Let $G$ be a finite abelian group contained in $\mathcal{C}$. Then $\alpha(G)=\frac{3}{4}$ and, by the fundamental theorem of finitely generated abelian groups, it follows that
$$G\cong\xmare{i=1}k G_i,$$
where $G_i$ is a finite abelian $p_i$-group of order $p_i^{n_i}$, with $n_i\in\mathbb{N}^*$, for all $i=1,2,\ldots, k$. Since we work with a family of finite groups $(G_i)_{i=\overline{1,k}}$ having coprime orders, we obtain
$$\alpha(G)=\prod\limits_{i=1}^k\alpha(G_i).$$
If there exists $i=1,2,\ldots,k$ such that $p_i$ is an odd prime number, then, by \textbf{Lemma 2.2}, we have $\alpha(G_i)\le\alpha(\mathbb{Z}_{p_i}^{n_i})=\frac{1+\frac{p_i^{n_i}-1}{p_i-1}}{p_i^{n_i}}$. It is easy to show that the function $f:[1,\infty)\longrightarrow \mathbb{R}, \ f(x)=\frac{1+\frac{p_i^x-1}{p_i-1}}{p_i^x}$ is strictly decreasing, so its maximum value is $\frac{2}{p_i}$. Therefore, $\alpha(G_i)\le\frac{2}{p_i}<\frac{3}{4}$. As a consequence, we obtain
$\frac{3}{4}=\alpha(G)=\prod\limits_{i=1}^k\alpha(G_i)<\frac{3}{4}$ and we arrive at a contradiction. Hence, $p_i=2$ for all $i=1,2,\ldots,k$. Then $G\cong \mathbb{Z}_{2^{d_1}}\times\mathbb{Z}_{2^{d_2}}\times\ldots\times \mathbb{Z}_{2^{d_k}}$, where $1\le d_1\le d_2\le\ldots\le d_k$.

We cannot have $d_k=1$ since this would imply that $G$ is an elementary abelian 2-group and, in this case, $\alpha(G)=1$, a contradiction. If $d_k\ge 3$, then $G$ would have a subgroup isomorphic to $\mathbb{Z}_{2^{d_k}}$ and we would obtain
$$\frac{3}{4}=\alpha(G)\le\alpha(\mathbb{Z}_{2^{d_k}})=\frac{d_k+1}{2^{d_k}}\le\frac{1}{2}<\frac{3}{4},$$ a contradiction. Note that in the above argument we used the fact that the function $g:[3,\infty)\longrightarrow \mathbb{R}$, given by $g(x)=\frac{x+1}{2^x}$, is strictly decreasing and its maximum value is $\frac{1}{2}$.

Hence, $d_k=2$, so we may assume that $d_1+d_2+\ldots+d_k=n$, where $n\ge 2$ is a positive integer. If $k=1$, we have $G\cong\mathbb{Z}_{2^{d_1}}$. Since $G\in\mathcal{C}$, we have
$$\alpha(G)=\frac{3}{4}\Longleftrightarrow \frac{d_1+1}{2^{d_1}}=\frac{3}{4}\Longleftrightarrow d_1=2,$$ which leads to $G\cong\mathbb{Z}_4$. If $1<k<n-1$, then $n>2$ and $G$ has a subgroup isomorphic to $\mathbb{Z}_4\times\mathbb{Z}_4$. We arrive at a contradiction as
$$\frac{3}{4}=\alpha(G)\le\alpha(\mathbb{Z}_4\times\mathbb{Z}_4)=\frac{5}{8}<\frac{3}{4}.$$
If $k=n-1$, then $n>2$ and it is obvious that $G\cong\mathbb{Z}_2^{n-2}\times \mathbb{Z}_4$.

Conversely, as we stated in the \textbf{Introduction}, for all $n\in\mathbb{N}$, we have
$$\alpha(\mathbb{Z}_2^n\times\mathbb{Z}_4)=\alpha(\mathbb{Z}_4)=\frac{3}{4}.$$ Therefore, it is true that $\mathbb{Z}_2^n\times\mathbb{Z}_4$, where $n\in\mathbb{N}$, are the only abelian groups contained in $\cal{C}$.
\hfill\rule{1,5mm}{1,5mm} \\

An alternative proof of \textbf{Theorem 2.3} may be given if one shows that $\alpha:\mathcal{A}\longrightarrow [0,1]$ is an injective function, where $\cal{A}$ is the class of finite abelian $p$-groups of a given order. In this sense, we recall that \textbf{Theorem 4.3} of \cite{12} provides an explicit formula for computing the number of cyclic subgroups of a given order that are contained in a finite abelian $p$-group. More exactly, for a finite abelian $p$-group $G\cong\mathbb{Z}_{p^{d_1}}\times\mathbb{Z}_{p^{d_2}}\times\ldots\times\mathbb{Z}_{p^{d_k}}$, where $1\le d_1\le d_2\le\ldots\le d_k$, the number of cyclic subgroups of order $p^i$, with $i=1,2,\ldots,d_k$, is
$$g_p^k(i)=\frac{p^ih_p^{k-1}(i)-p^{i-1}h_p^{k-1}(i-1)}{p^i-p^{i-1}},$$
where
$$h_p^{k-1}(i)=\begin{cases} p^{(k-1)i}  &\mbox{, \ } 0\le i\le d_1 \\  p^{(k-2)i+d_1} &\mbox{, \ } d_1\le i\le d_2 \\ \vdots \\ p^{d_1+d_2+\ldots +d_{k-1}} &\mbox{, \ } d_{k-1}\le i\end{cases}.$$
Taking $d_0=0$, adding the quantities $g_p^k(i)$, where $i=1,2,\ldots,d_k$, and counting the trivial subgroup, we obtain the following result.\\

\textbf{Theorem 2.4.} \textit{Let $G\cong\mathbb{Z}_{p^{d_1}}\times\mathbb{Z}_{p^{d_2}}\times\ldots\times\mathbb{Z}_{p^{d_k}}$, where $1\le d_1\le d_2\le\ldots\le d_k$ and $p$ is a prime number. Then the total number of cyclic subgroups of $G$ is
$$1+\frac{1}{p-1}\sum\limits_{i=0}^{k-2}p^{d_0+d_1+\ldots+d_i}\frac{p^{k-i}-1}{p^{k-i-1}-1}(p^{(k-i-1)d_{i+1}}-p^{(k-i-1)d_i})+(d_k-d_{k-1})p^{d_0+d_1+\ldots+d_{k-1}}.$$ }

Besides the fact that \textbf{Theorem 2.4} provides an explicit formula that allow us to compute the total number of cyclic subgroups of any finite abelian $p$-group, we consider that this result may be useful for proving the following conjecture.\\

\textbf{Conjecture 2.5.} \textit{The function $\alpha:\mathcal{A}\longrightarrow [0,1]$ is injective.}\\

We are ready to prove our main result which describes the characteristics of the groups that belong to $\cal{C}$.\\

\textbf{Theorem 2.6.} \textit{Let $G\in\cal{C}$. Then $G$ is a 2-group with $G'=\Phi(G)$ or there is $n\in\mathbb{N}$ such that $\frac{G}{G'}\cong \mathbb{Z}_2^n\times \mathbb{Z}_4$ and $G'$ is elementary abelian.}

\textbf{Proof.} Let $G\in\cal{C}$. Then $G$ is a finite nilpotent group, so it is isomorphic to the direct product of its Sylow subgroups. Since $\alpha(G)=\frac{3}{4}$ and the Sylow subgroups are having coprime orders, following a similar argument with the one that was used in the proof of \textbf{Theorem 2.3}, one obtains that $G$ is a 2-group. Also, remark that we have
$$\frac{3}{4}=\alpha(G)\le\alpha(\frac{G}{G'}),$$ so we distinguish the following two cases.

If $\alpha(\frac{G}{G'})>\frac{3}{4}$, then $\frac{G}{G'}$ is isomorphic to one of the groups described by \textbf{Theorem 5} of \cite{4}. Since $\frac{G}{G'}$ is an abelian group, the only possible choice from the respective classification is $\frac{G}{G'}\cong \mathbb{Z}_2^n$, where $n$ is a positive integer. It follows that $\Phi(G)\subseteq G'$. The converse inclusion is well known, so the equality $G'=\Phi(G)$ holds.

If $\alpha(\frac{G}{G'})=\frac{3}{4}$, then, according to \textbf{Theorem 2.3}, there exists $n\in\mathbb{N}$ such that $\frac{G}{G'}\cong \mathbb{Z}_2^n\times\mathbb{Z}_4$. Also, in this case, we have $\alpha(G)=\alpha(\frac{G}{G'})$. The proof is complete since the last equality implies that $G'$ is an elementary abelian 2-group.
\hfill\rule{1,5mm}{1,5mm} \\

We end this section by stating that our last result may be improved if one would classify the finite groups having the properties indicated by \textbf{Theorem 2.6}. This may lead to a complete determination of the class $\cal{C}$.

\section{Connections between $\cal{C}$ and other classes of finite groups}

Our next aim is to study if there are some well known types of finite groups which belong to $\cal{C}$. We remark that we can limit our study to finite 2-groups since this is one of the features of the groups that are contained in $\cal{C}$. Hence, we will focus on the following classes of finite groups: (almost) extraspecial 2-groups, generalized dicyclic 2-groups, generalized dihedral 2-groups and 2-groups possessing a cyclic maximal subgroup.

For a finite group $G$, we denote by $I(G)$ the number of involutions of $G$, i.e. the number of elements of order 2 of $G$. Our next result shows that, under certain assumptions, there is a connection between the class $\cal{C}$ and $I(G)$.\\

\textbf{Proposition 3.1.} \textit{Let $n\ge 2$ be a positive integer and let $G$ be a 2-group of order $2^n$ having $exp(G)=4$. Then $G\in\cal{C}$ if and only if $I(G)=2^{n-1}-1$.}

\textbf{Proof.} Recall that we denoted by $n_{2^i}(G)$ the number of cyclic subgroups of order $2^i$ of $G$, for $i=1,2$. We remark that $I(G)=n_2(G)$. Also, we have
$$2^n=1+n_2(G)+2n_4(G).$$
Therefore,
\begin{align*}
G\in\mathcal{C} & \Longleftrightarrow \alpha(G)=\frac{3}{4} \\ &\Longleftrightarrow \frac{1+n_2(G)+n_4(G)}{2^n}=\frac{3}{4} \\ &
 \Longleftrightarrow \frac{2^n-n_4(G)}{2^n}=\frac{3}{4} \\ & \Longleftrightarrow n_4(G)=2^{n-2}\\ &\Longleftrightarrow n_2(G)=2^{n-1}-1,
\end{align*}
as desired.
\hfill\rule{1,5mm}{1,5mm} \\

\textbf{Proposition 3.1} indicates that the $2$-groups of exponent 4 that belong to $\cal{C}$ are exactly the $2$-groups of exponent 4 described in \cite{7}. Also, this result characterizes the appartenance to $\cal{C}$ of all classes of finite 2-groups of exponent 4. Two of these classes are formed by extraspecial 2-groups and almost extraspecial 2-groups, respectively. Before we study the connection between $\cal{C}$ and the previously mentioned classes of groups, we recall some theoretical aspects related to central products of groups.\\

Let $G$ be a finite group and let $H_1,H_2$ be two of its subgroups. Then $G$ is the internal central product of $H_1$ and $H_2$ if $G=H_1H_2$ and $[H_1,H_2]=\lbrace 1\rbrace$. We denote this fact by $G=H_1*H_2$. Also, for a positive integer $r$, we denote by $H_1^{*r}$ the central product of $r$ copies of $H_1$. A connection with the usual external direct product of the subgroups $H_1$ and $H_2$ is given by \textbf{Theorem 3.4} of \cite{6}. In some words, this result shows how one can obtain the external central product, which is isomorphic to $G=H_1*H_2$, by quotiening the external direct product $H_1\times H_2$ by a normal subgroup $Z$. Hence, $G=H_1*H_2\cong \frac{H_1\times H_2}{Z}$. Also, \textbf{Example 3.5} of the same paper shows how one may obtain the structure of the subgroup $Z$ starting with the internal central product $D_8*\mathbb{Z}_4$, but the same process can be applied in more general situations.

We recall that a 2-group $G$ is called
\begin{itemize}
\item[--] extraspecial if $G'=\Phi(G)=Z(G)\cong\mathbb{Z}_2$;
\item[--] almost extraspecial if $G'=\Phi(G)\cong\mathbb{Z}_2$ and $Z(G)\cong\mathbb{Z}_4.$
\end{itemize}
Moreover, according to \textbf{Theorem 2.3} of \cite{1},
\begin{itemize}
\item[--] if $G$ is an extraspecial 2-group, then there is a positive integer $r$ such that $|G|=2^{2r+1}$ and $G\cong D_8^{*r}$ or $G\cong Q_8*D_8^{*(r-1)}$;
\item[--] if $G$ is an almost extraspecial 2-group, then there is a positive integer $r$ such that $|G|=2^{2r+2}$ and $G\cong D_8^{*r}*\mathbb{Z}_4$.
\end{itemize}
We remark that excepting $Q_8$, any (almost) extraspecial 2-group $G$ can be written as a central product $D_8*G_1$. More exactly, for a positive integer $r$,
\begin{itemize}
\item[--] if $G$ is an extraspecial 2-group, then $G_1\cong D_8^{*(r-1)}$ or $G_1\cong D_8^{*(r-2)}*Q_8$;
\item[--] if $G$ is an almost extraspecial 2-group, then $G_1\cong D_8^{*(r-1)}*\mathbb{Z}_4$.
\end{itemize}
Let $D_8=\langle x, y \ | \ x^4=y^2=1, yx=x^3y\rangle$. If $G$ is an extraspecial 2-group, then starting with $G_1\cong D_8$ or $G_1\cong Q_8$, passing to the external central product $\frac{D_8\times G_1}{Z}$ and increasing the number of central factors of type $D_8$ after each such step, we infer that $exp(G)=4$, $Z(G_1)=\langle a\rangle\cong\mathbb{Z}_2$ and $Z=\langle (x^2,a)\rangle\cong\mathbb{Z}_2$. If $G$ is an almost extraspecial 2-group, then starting with $G_1\cong \mathbb{Z}_4$, passing to the external central product $\frac{D_8\times G_1}{Z}$ and increasing the number of central factors of type $D_8$, it follows that $exp(G)=4$, $Z(G_1)=\langle a\rangle\cong\mathbb{Z}_4$ and $Z=\langle (x^2,a^2)\rangle\cong\mathbb{Z}_2$. We note that, for ease in writing, we used the same letter to denote the generator of $Z(G_1)$ in both cases. Our next result is relevant for finding the (almost) extraspecial 2-groups contained in $\cal{C}$.\\

\textbf{Lemma 3.2.} \textit{Let $n\ge 4$ be a positive integer and let $G\cong D_8*G_1$ be a finite (almost) extraspecial 2-group of order $2^n$. Then $n_2(G)=2^{n-2}+2n_2(G_1)+1.$}

\textbf{Proof.} Let $G\cong D_8*G_1$ be a finite (almost) extraspecial 2-group or order $2^n$, where $n\ge 4$ is a positive integer. Then, $G\cong\frac{D_8\times G_1}{Z}$ and:
\begin{itemize}
\item[--] $Z(G_1)=\langle a\rangle\cong\mathbb{Z}_2$ and $Z=\langle (x^2,a)\rangle\cong\mathbb{Z}_2$, if $G$ is an extraspecial 2-group;
\item[--] $Z(G_1)=\langle a\rangle\cong\mathbb{Z}_4$ and $Z=\langle (x^2,a^2)\rangle\cong\mathbb{Z}_2$, if $G$ is an almost extraspecial 2-group.
\end{itemize}
To find the total number of cyclic subgroups of order $2$ of $G$, i.e. the quantity $n_2(G)$, the first step is to count the total number of elements $(u,v)\in D_8\times G_1$ such that $(u,v)Z$ has order 2 in $G$. Hence, for an arbitrary element $(u,v)\in D_8\times G_1$, we have
$$ord((u,v)Z)=2 \Longleftrightarrow \left\{\begin{array}{ll}
(u,v)\not\in Z  \\
(u^2,v^2) \in Z
\end{array} \right.\Longleftrightarrow \left\{\begin{array}{ll}
(u,v)\ne (1,1), (x^2,a^i)  \\
(u^2,v^2)=(1,1) \text{ \ or \ } (u^2,v^2)=(x^2,a^i)
\end{array} \right.,$$
where $i=1$, if $G$ is an extraspecial 2-group, and $i=2$, if $G$ is an almost extraspecial 2-group. There are $6(1+n_2(G_1))$ solutions for the equation $(u^2,v^2)=(1,1)$ since there are 6 elements $u\in D_8$ such that $u^2=1$, and $1+n_2(G_1)$ elements $v\in G_1$ satisfying $v^2=1$. Similarly, for $i=1,2$, the equation $(u^2,v^2)=(x^2,a^i)$ has $2(2^{n-2}-1-n_2(G_1))$ solutions. Since $(u,v)\ne (1,1), (x^2,a^i)$ and these pairs are solutions for the equation $(u^2,v^2)=(1,1)$, there are
$$6(1+n_2(G_1))+2(2^{n-2}-1-n_2(G))-2=2^{n-1}+4n_2(G_1)+2$$
elements $(u,v)\in D_8\times G_1$ such that $ord((u,v)Z)=2$. Since $G\cong D_8*G_1\cong\frac{D_8\times G_1}{Z}$ and $|Z|=2$, the number of cyclic subgroups of order $2$ of $G$ is
$$n_2(G)=\frac{2^{n-1}+4n_2(G_1)+2}{2}=2^{n-2}+2n_2(G)+1.$$
\hfill\rule{1,5mm}{1,5mm} \\

Under the hypotheses of \textbf{Lemma 3.2}, since $n_2(G)$ is computed, one can easily obtain that
$$n_4(D_8*G_1)=3\cdot 2^{n-3}-n_2(G_1)-1.$$
Using this explicit result, it follows that
$$n_4(D_8*Q_8)=10 \text { \ and \ } n_4(D_8*D_8)=6.$$
The same numbers were obtained in the proof of \textbf{Proposition 3.13} of \cite{3}. Also, we have
$$n_2(D_8*\mathbb{Z}_4)= 7 \text{ \ and \ } n_4(D_8*\mathbb{Z}_4)=4,$$
the same results being indicated by \textbf{Example 4.5} of \cite{6}.

We are ready to study the connection between $\cal{C}$ and the class of finite (almost) extraspecial 2-groups.\\

\textbf{Theorem 3.3.} \textit{a) There are no finite extraspecial 2-groups contained in $\cal{C}.$}\\
\hspace*{2.9cm} \textit{b) Any finite almost extraspecial 2-group belongs to $\cal{C}.$}

\textbf{Proof.} Let $G$ be a finite (almost) extraspecial 2-group of order $2^n$, where $n\ge 4$ is a positive integer. We recall that $exp(G)=4$, so we can use \textbf{Proposition 3.1} in both cases.

a) If $G$ is extraspecial, then $G\cong D_8^{*r}$ or $G\cong Q_8*D_8^{*(r-1)}$, where $r$ is a positive integer. Since $Q_8$ is not contained in $\cal{C}$, it is sufficient to study the appartenance to $\cal{C}$ of $G\cong D_8*G_1$, where $G_1\cong D_8^{*(r-1)}$ or $G_1\cong D_8^{*(r-2)}*Q_8$. Using \textbf{Proposition 3.1}, \textbf{Lemma 3.2} and the fact that $|G_1|=2^{n-2}$, we have
\begin{align*}
G\cong D_8*G_1\in \mathcal{C}& \Longleftrightarrow n_2(D_8*G_1)=2^{n-1}-1 \\ &\Longleftrightarrow 2^{n-2}+2n_2(G_1)+1=2^{n-1}-1 \\ &\Longleftrightarrow n_2(G_1)=2^{n-3}-1 \\ & \Longleftrightarrow G_1\in\cal{C}.
\end{align*}
Repeating this argument for a finite number of times, we obtain that
$$G\cong D_8*G_1\in\mathcal{C} \Longleftrightarrow D_8\in\mathcal{C} \text{ \ or \ } Q_8\in\mathcal{C}.$$
Since $D_8$ and $Q_8$ do not belong to $\cal{C}$, it follows that $G\not\in\cal{C}.$

b) If $G$ is almost extraspecial, then $G\cong D_8^{*r}*\mathbb{Z}_4$, where $r$ is a positive integer. In this case, $G\cong D_8*G_1$, where $G_1\cong D_8^{*(r-1)}*\mathbb{Z}_4$. Using the same reasoning as in the proof of a), we have
$$G\cong D_8*G_1\in\mathcal{C} \Longleftrightarrow G_1\in \mathcal{C}.$$
Repeating the same steps for a finite number of times, we are led to
$$G\cong D_8*G_1\in\mathcal{C}\Longleftrightarrow \mathbb{Z}_4\in\mathcal{C}.$$
Therefore, since the right statement of the above equivalence is true, all finite almost extraspecial 2-groups are contained in $\cal{C}.$
\hfill\rule{1,5mm}{1,5mm} \\

The next step is to find the generalized dicyclic 2-groups that are contained in $\cal{C}$. Let $n\ge 2$ be a positive integer and $A$ be an abelian group of order $2^{n-1}$. Then a generalized dicyclic group of order $2^n$ has the following structure
$$Dic_{2^n}(A)=\langle A, \gamma \ | \ \gamma^4=1, \gamma^2\in A\setminus\lbrace 1\rbrace, \gamma g=g^{-1}\gamma, \ \forall \ g\in A\rangle.$$
For more details about this class of groups, we refer the reader to \cite{14}, where some probabilistic aspects associated to (generalized) dicyclic groups were studied.

We remark that if $A\cong \mathbb{Z}_2^{n-1}$, then $\gamma g=g\gamma, \ \forall \ g\in A$. This further implies that $Dic_{2^n}(A)$ is a finite abelian 2-group, so all its subgroups are normal. If we denote the $n-1$ generators of $A$ with $a_i$, for $i=1,2,\ldots,n-1$, without loss of generality, we may choose $\gamma^2$ to be $a_{n-1}$. Then the subgroups $H=\langle a_1\rangle\times\langle a_2\rangle\times\ldots\times\langle a_{n-2}\rangle\cong \mathbb{Z}_2^{n-2}$ and $K=\langle\gamma\rangle\cong\mathbb{Z}_4$ of $Dic_{2^n}(A)$ have trivial intersection and $Dic_{2^n}(A)=HK$. Hence, $Dic_{2^n}(A)=H\times K\cong \mathbb{Z}_2^{n-2}\times \mathbb{Z}_4.$ Having in mind this isomorphism, we are able to prove the following result which states that the generalized dicyclic 2-groups that belong to $\cal{C}$ are isomorphic to the abelian groups contained in $\cal{C}$.\\

\textbf{Theorem 3.4} \textit{The only generalized dicyclic 2-groups contained in $\mathcal{C}$ are isomorphic to $\mathbb{Z}_2^n\times\mathbb{Z}_4,$ where $n\in\mathbb{N}$.}

\textbf{Proof.} Let $A$ be an abelian group of order $2^{n-1}$, where $n\ge 2$ is a positive integer. We denote the non-trivial elements of $A$ by $g_i$, for $i=1,2,\ldots,2^{n-1}-1$. Then
$$Dic_{2^n}(A)=\lbrace 1,g_1,g_2,\ldots, g_{2^{n-1}-1},\gamma, g_1\gamma,g_2\gamma,\ldots,g_{2^{n-1}-1}\gamma\rbrace.$$
It is clear that $L_1(A)\subset L_1(Dic_{2^n}(A))$. Also, since $ord(g_i\gamma)=4, \ \forall \ i=1,2,\ldots,2^{n-1}-1$ and $\langle g_i\gamma\rangle\cap\langle g_j\gamma\rangle=\langle\gamma^2\rangle$, $\forall \ i,j\in\lbrace 1,2,\ldots,2^{n-1}-1\rbrace$ with $i\ne j$, $Dic_{2^n}(A)$ has other $2^{n-2}$ cyclic subgroups of order 4. If we assume that $Dic_{2^n}(A)$ possesses another cyclic subgroup $H$ besides the ones that were already indicated, then there exists $i=1,2,\ldots,2^{n-1}-1$ such that $H\cap\langle g_i\gamma\rangle=\langle g_i\gamma\rangle$. Then, $H$ should be generated by an element $g$ of $Dic_{2^n}(A)$ with $ord(g)\ge 8$. But such elements may be contained only in $A$ and this implies that $\langle g_i\gamma\rangle= H\cap \langle g_i\gamma\rangle\le \langle\gamma^2\rangle$, a contradiction. Therefore, we have
$$|L_1(Dic_{2^n}(A))|=|L_1(A)|+2^{n-2}.$$

Assume that $Dic_{2^n}(A)\in\cal{C}$.
It follows that
$$\alpha(Dic_{2^n}(A))=\frac{3}{4}\Longleftrightarrow \frac{|L_1(A)|+2^{n-2}}{2^n}=\frac{3}{4}\Longleftrightarrow |L_1(A)|=2^{n-1}.$$
Hence, we must determine the finite abelian groups $A$ satisfying $|A|=|L_1(A)|=2^{n-1}$. The following set of conditions
$$\left\{\begin{array}{ll}
2^{n-1}=1+n_2(G)+2n_4(G)+\ldots +2^{m-1}n_{2^m}(G)\\
2^{n-1}=1+n_2(G)+n_4(G)+\ldots +n_{2^m}(G)
\end{array} \right.$$
holds, where $exp(A)=2^m$. Since $n\ge 2$, we have $m\ge1$. If we assume that $m\ge 2$, then the above set of conditions leads to $n_4=n_8=\ldots=n_{2^m}=0$, a contradiction. Therefore, $A$ is an abelian 2-group of order $2^{n-1}$ and its exponent is 2. Then $A\cong\mathbb{Z}_2^{n-1}$, so, according to our remarks that were made before we started this proof, we have $Dic_{2^n}(A)\cong\mathbb{Z}_2^{n-2}\times\mathbb{Z}_4$.

For the converse, we already indicated that for all $n\in\mathbb{N}$, the abelian 2-groups $\mathbb{Z}_2^n\times\mathbb{Z}_4$ belong to $\cal{C}$. Hence, our proof is complete.
\hfill\rule{1,5mm}{1,5mm} \\

Our following aim is to establish which are the generalized dihedral 2-groups that are contained in $\cal{C}$. We denote by $y$ the generator of the cyclic group $\mathbb{Z}_2$. We begin by recalling that for an abelian group $G$, the generalized dihedral group associated to $G$ is $D(G)=G\rtimes_{\varphi}\mathbb{Z}_2$, where $\varphi:\mathbb{Z}_2\longrightarrow Aut(G)$ is a homomorphism given by:
$$\left\{\begin{array}{ll}
\varphi(1)=1_{G}\\
\varphi(y)=\varphi_{y}, \ \varphi_{y}(g)=g^{-1}, \ \forall \ g\in G
\end{array} \right..$$
A presentation of the generalized dihedral group $D(G)$ is the following one
$$D(G)=\langle G,y \ | \ y^2=1, ygy=g^{-1}, \ \forall \ g\in G\rangle.$$
It is known that $D(\mathbb{Z}_n)=D_{2n}$ for any positive integer $n\ge2$. Other properties of generalized dihedral groups and the dihedralization of several finite abelian groups are presented in \cite{2}. To find the generalized dihedral 2-groups that are contained in $\cal{C}$, we will need the following preliminary result which provides a classification of the abelian 2-groups having $\alpha(G)=\frac{1}{2}.$\\

\textbf{Lemma 3.5.} \textit{Let $G$ be a finite abelian 2-group. Then $\alpha(G)=\frac{1}{2}$ if and only if there is $n\in\mathbb{N}$ such that $G\cong \mathbb{Z}_2^n\times\mathbb{Z}_8.$}

\textbf{Proof.} Let $m$ be a positive integer and $G$ be a finite abelian 2-group of order $2^m$ such that $\alpha(G)=\frac{1}{2}$. Firstly, we show that
$$exp(G)\le 8.$$
Indeed, if $G\cong\mathbb{Z}_{2^{d_1}}\times\mathbb{Z}_{2^{d_2}}\times\ldots\times\mathbb{Z}_{2^{d_k}}$ where $1\le d_1\le d_2\le\ldots\le d_k$, then
$$\frac{d_k+1}{2^{d_k}}=\alpha(\mathbb{Z}_{2^{d_k}})=\alpha\bigg(\frac{G}{\mathbb{Z}_{2^{d_1}}\times\mathbb{Z}_{2^{d_1}}\times\ldots\times \mathbb{Z}_{2^{d_{k-1}}}}\bigg)\ge\alpha(G)=\frac{1}{2}.$$ This leads to $d_k\le 3$, and, consequently, $exp(G)\le 8$.

Obviously, we cannot have $exp(G)=2$ since this would imply that $G\cong\mathbb{Z}_2^k$ and $\alpha(G)=1$, a contradiction. Assume that $exp(G)=4$. Then, the equality $\alpha(G)=\frac{1}{2}$ leads to $|L_1(G)|=2^{m-1}$. Hence, the following set of conditions holds:
$$\left\{\begin{array}{ll}
2^m=1+n_2(G)+2n_4(G)\\
2^{m-1}=1+n_2(G)+n_4(G)
\end{array} \right..$$
We infer that $n_4(G)=2^{m-1}$ and, consequently, we obtain $1+n_2(G)=0$, a contradiction. Therefore, the exponent of $G$ is 8, so $G\cong \mathbb{Z}_2^n\times\mathbb{Z}_4^a\times\mathbb{Z}_8^b$, where $a,b,n\in\mathbb{N}$. Note that at least one of $a$ and $b$ is strictly positive, since, otherwise, we would have $\alpha(G)=1$, a contradiction.  Also, the following set of conditions holds:
$$\left\{\begin{array}{ll}
2^m=1+n_2(G)+2n_4(G)+4n_8(G)\\
2^{m-1}=1+n_2(G)+n_4(G)+n_8(G)
\end{array} \right..$$
The above equalities imply that $n_4(G)+3n_8(G)=2^{m-1}.$ This further leads to
$$1+n_2(G)=2n_8(G).$$
Using \textbf{Theorem 4.3} of \cite{12}, the numbers of cyclic subgroups of order 2 and 8, respectively, are
$$n_2(G)=2^{n+a+b}-1 \text{ \ and \ } n_8(G)= 2^{n+2a+2b-2}(2^b-1).$$
Then, we obtain
\begin{align*}
1+n_2(G)=2n_8(G) &\Longleftrightarrow 2^{n+a+b}=2^{n+2a+2b-1}(2^b-1)\\
&\Longleftrightarrow 1=2^{a+b-1}(2^b-1)\\
&\Longleftrightarrow a=0 \text{ \ and \ } b=1.
\end{align*}
Hence, $G\cong \mathbb{Z}_2^n\times\mathbb{Z}_8$, as desired.

Conversely, if $G\cong\mathbb{Z}_2^n\times\mathbb{Z}_8$, where $n\in\mathbb{N}$, we have
$$\alpha(G)=\alpha(\mathbb{Z}_2^n\times\mathbb{Z}_8)=\alpha(\mathbb{Z}_8)=\frac{1}{2},$$
and our proof is complete.
\hfill\rule{1,5mm}{1,5mm} \\

Our next result examines the connection between $\cal{C}$ and the class of finite generalized dihedral 2-groups.\\

\textbf{Theorem 3.6.} \textit{The only finite generalized dihedral 2-groups that belong to $\cal{C}$ are isomorphic to $\mathbb{Z}_2^n\times D_{16}$, where $n\in\mathbb{N}$.}

\textbf{Proof.} Let $G$ be a finite abelian group such that $D(G)$ is a generalized dihedral group contained in $\cal{C}$. It is known that for any subgroup $H$ of $G$, the generalized dihedral group $D(G)$ has one subgroup isomorphic to $H$ and $[G:H]$ subgroups isomorphic to $D(H)$. However, we are interested in counting only the cyclic subgroups of $D(G)$. It is clear that
$$L_1(G) \cup \lbrace \langle gy\rangle \ | \ g\in G\rbrace\subseteq L_1(D(G)).$$
Assume that $D(G)$ has another cyclic subgroup $H$. Then $H$ contains at least one subgroup of type $\langle gy\rangle$, where $g\in G$, and $|H|\ge 4$. This implies that $H$ must be generated by an element $h\in G$ and we would have $$\langle gy\rangle=H\cap\langle gy\rangle=\langle h\rangle\cap\langle gy\rangle\subseteq G\cap\langle gy\rangle=\lbrace 1\rbrace,$$
a contradiction. Then,
$$|L_1(D(G))|=|L_1(G)|+|G|.$$
Since $D(G)$ is a finite 2-group, it follows that $G$ is a finite abelian 2-group. Moreover,
$$D(G)\in\mathcal{C}\Longleftrightarrow \alpha(D(G))=\frac{3}{4}\Longleftrightarrow \frac{|L_1(G)|+|G|}{2|G|}=\frac{3}{4}\Longleftrightarrow \alpha(G)=\frac{1}{2}.$$
According to \textbf{Lemma 3.5}, there is $n\in\mathbb{N}$ such that $G\cong\mathbb{Z}_2^n\times\mathbb{Z}_8.$ Using \textbf{Theorem 5.1} of \cite{2}, we obtain
$$D(G)\cong D(\mathbb{Z}_2^n\times\mathbb{Z}_8)\cong\mathbb{Z}_2^n\times D(\mathbb{Z}_8)\cong\mathbb{Z}_2^n\times D_{16}.$$

Conversely, we have
$$\alpha(\mathbb{Z}_2^n\times D_{16})=\alpha(D_{16})=\frac{3}{4},$$
a fact which completes our proof.
\hfill\rule{1,5mm}{1,5mm} \\

The last connection that we study is between $\cal{C}$ and the class of finite 2-groups possessing a cyclic maximal subgroup. Let $n\ge 3$ be a positive integer. The class of finite 2-groups possessing a maximal subgroup which is cyclic contains abelian groups of type $\mathbb{Z}_2\times\mathbb{Z}_{2^{n-1}}$ which, according to \textbf{Theorem 2.3}, belong to $\cal{C}$ if and only if $n=3$. Hence, it is sufficient to study the appartenance of the non-abelian 2-groups possessing a cyclic maximal subgroup to $\cal{C}$. \textbf{Theorem 4.1} of \cite{9}, II, provides a complete classification of the non-abelian 2-groups containing a maximal subgroup which is cyclic. They are isomorphic to
\begin{itemize}
\item[--] the modular 2-group $$M(2^n)=\langle x,y \ | \ x^{2^{n-1}}=y^2=1, y^{-1}xy=x^{2^{n-2}+1}\rangle, n\geq 4,$$
\item[--] the dihedral group $$D_{2^n}=\langle x,y \ | \ x^{2^{n-1}}=y^2=1, yxy=x^{-1}\rangle, n\ge 3,$$
\item[--] the generalized quaternion group $$Q_{2^n}=\langle x,y \ | \ x^{2^{n-1}}=y^4=1, yxy^{-1}=x^{-1}\rangle, n\geq 3,$$
\item[--] the quasi-dihedral group $$S_{2^n}=\langle x,y \ | \ x^{2^{n-1}}=y^2=1, yxy=x^{2^{n-2}-1}\rangle, n\geq 4.$$
\end{itemize}
In \cite{13}, the number of cyclic subgroups of each of these groups is indicated. More exactly, we have
$$|L_1(M(2^n))|=2n, \ \ |L_1(D_{2^n})|=2^{n-1}+n, \ \ |L_1(Q_{2^n})|=2^{n-2}+n, \ \ |L_1(S_{2^n})|=3\cdot 2^{n-3}+n.$$
Solving the equation $\alpha(G)=\frac{3}{4}$, where $G$ is isomorphic to one of the above groups, over the positive integers, we obtain
$$\alpha(M(2^n))=\frac{3}{4} \Longleftrightarrow \frac{2n}{2^n}=\frac{3}{4}\Longleftrightarrow n=3,$$
$$\alpha(D_{2^n})=\frac{3}{4} \Longleftrightarrow \frac{2^{n-1}+n}{2^n}=\frac{3}{4}\Longleftrightarrow n=4,$$
$$\alpha(Q_{2^n})=\frac{3}{4} \Longleftrightarrow \frac{2^{n-2}+n}{2^n}=\frac{3}{4}\Longleftrightarrow n=1 \text{ \ or \ } n=2,$$
$$\alpha(S_{2^n})=\frac{3}{4} \Longleftrightarrow \frac{3\cdot 2^{n-3}+n}{2^n}=\frac{3}{4}\Longleftrightarrow n=3.$$
This leads us to our final result of this paper.\\

\textbf{Theorem 3.7.} \textit{The dihedral group $D_{16}$ is the only non-abelian 2-group that possesses a cyclic maximal subgroup and belongs to $\cal{C}$.\\}

We end our paper by noticing that, up to a direct factor of type $\mathbb{Z}_2^n,$ where $n\ge 1$ is a positive integer, $D_{16}$ is:
\begin{itemize}
\item[--] the only finite generalized dihedral 2-group contained in $\cal{C}$;
\item[--] the only finite non-abelian 2-group that possesses a cyclic maximal subgroup and belongs to $\cal{C}$.
\end{itemize}

\section{Further research}

We studied the connections between $\cal{C}$ and other classes of finite groups and we indicated the characteristics of the groups that belong to $\cal{C}$. However, we did not manage to provide a complete classification of the groups contained in $\cal{C}$. Hence, besides \textbf{Conjecture 2.5}, we indicate the following open problems:\\

\textbf{Problem 4.1.} Study the connections between other well known classes of finite 2-groups and $\cal{C}$.\\

\textbf{Problem 4.2.} Classify all finite groups that are contained in $\cal{C}$.\\

\textbf{Problem 4.3.} Classify all finite groups satisfying $\alpha(G)=\frac{3}{4}.$\\

\textbf{Problem 4.4.} Study the density of the set $\lbrace \alpha(G) \ | \ G=\text{finite group}\rbrace$ in $[0,\frac{3}{4}]$.

\vspace*{3ex}
\small

\begin{minipage}[t]{7cm}
Marius T\u arn\u auceanu \\
Faculty of  Mathematics \\
"Al.I. Cuza" University \\
Ia\c si, Romania \\
e-mail: {\tt tarnauc@uaic.ro}
\end{minipage}
\hfill
\begin{minipage}[t]{7cm}
Mihai-Silviu Lazorec \\
Faculty of  Mathematics \\
"Al.I. Cuza" University \\
Ia\c si, Romania \\
e-mail: {\tt mihai.lazorec@student.uaic.ro}
\end{minipage}
\end{document}